\documentclass[a4paper,12pt]{amsart}

\usepackage{amssymb}
\usepackage{amsmath}
\usepackage{graphicx}

\title{Escher degree of non-periodic L-tilings by 2 prototiles.}

\author{Kazushi Ahara, Mami Murata, and Anno Ojiri}
\address{Department of Mathematics, Meiji University, 1-1-1 Higashi-Mita, Tama-ku, Kawasaki, Kanagawa, 214-8571, Japan}
\email{kazuaha63@hotmail.co.jp}


\newtheorem{thm}{Theorem}[section]

\newtheorem{lem}[thm]{Lemma}

\theoremstyle{definition}
\newtheorem{definition}[thm]{Definition}
\newtheorem{example}[thm]{Example}
\newtheorem{remark}[thm]{Remark}

\begin{document}

\begin{abstract}
For a given tiling of the euclidean plane ${\bf E}^2$, we call the degree of freedom of perturbed edges of prototiles {\it escher degree}.  In this paper we consider non-periodic L-tilings by $2$ prototiles and obtain the escher degree of them.
\end{abstract}

\maketitle


\section{introduction}

A non-periodic L-tiling is a limit of the sequence of tilings as shown in Figure 1.  (Often this is called {\it a chair tiling}.)  
It is well known that this is a tiling of the euclidean plane ${\bf E}^2$ and that it has no periodicity of parallel transformation.
\begin{figure}
\includegraphics[]{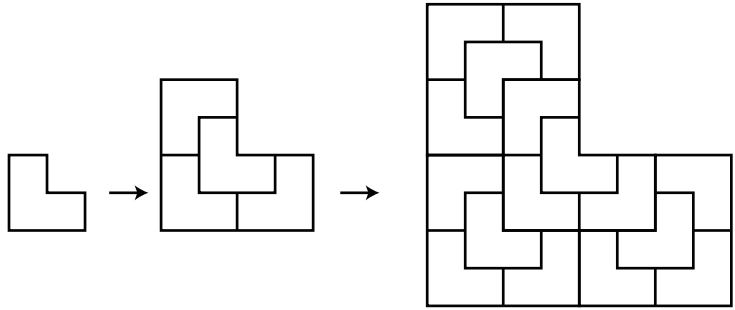}
\caption{non-periodic L-tiling}
\end{figure}

Sugihara \cite{Sugihara} introduces escherization of a plane tiling. Let $\mathcal T$ be a tiling of ${\bf E}^2$. 
If we have a finite set $\mathcal{S}=\{\alpha_1,\alpha_2,\cdots,\alpha_\ell\}$ of connected regions, and each tile of $\mathcal T$ is (orientation preserving) congruent to one of $\alpha_1,\alpha_2,\cdots,\alpha_\ell$, then we call $\mathcal{S}$ {\it the tile set} and $\alpha_1,\alpha_2,\cdots,\alpha_\ell$ {\it prototiles}.
If we perturb some of edges of prototiles and get another tiling of the plane, we call the process of perturbation {\it escherization} of the tiling $\mathcal T$.  This is a famous technique in artworks of M. C. Escher.
%

For example,  see Figure 2.  The left figure is a tiling by one parallelogram. We can perturb the horizontal edges and slanted edges independently as in the right figure.
\begin{figure}
\includegraphics[]{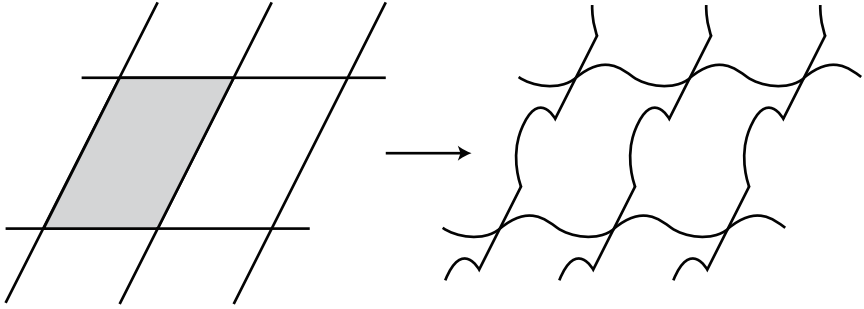}
\caption{perturbation of edges of parallelogram}
\end{figure}

In this paper, we determine the escher degree of L-tilings, that is, the degree of freedom of perturbed edges of prototiles of L-tilings.
If the tile set of an L-tiling consists of one prototile, then the escher degree is one.  This is shown in Theorem A.  If the tile set of an L-tiling consists of two prototiles, then we show that there are 6 types of non-trivial tilings. (A non-trivial tiling is a tiling whose escher degree is more than $1$.)  This is shown in Theorem B and Theorem C.
%

We have new kinds of tile sets in these theorems, but none of them are aperiodic.  This is shown in Appendix (2).
%

This paper is organized as follows.  In Section 2 we prepare some notations and basic lemmas.  In Section 3 we consider an L-tiling by one prototile.  In Section 4,5 we consider L-tilings by two prototiles.  
In appendix, we show figures of tilings. 
%


\section{preliminary}

In this section, first we introduce a non-periodic hierarchical tiling.  Let $\mathcal S=\{\alpha_1,\alpha_2,\cdots,\alpha_\ell\}$ be a set of connected regions and $\lambda>1$ a constant.  Let $\alpha'_i$ ($i=1,2,\cdots,\ell$) be a $\lambda$ scale-up copy of $\alpha_i$.
%

Suppose that each $\alpha'_i$ ($i=1,2,\cdots ,\ell$) can be tiled by prototiles $\alpha_1,\cdots , \alpha_\ell$.  That is, each $\alpha'_i$ can be divided into some of copies of  $\alpha_1,\cdots , \alpha_\ell$.  
Let $\alpha''_i$ be a $\lambda$ scale-up copy of $\alpha'_i$.  Then in the same way $\alpha''_i$ can be tiled by $\alpha'_1,\cdots , \alpha'_\ell$.  Substituting $\alpha_1,\cdots , \alpha_\ell$ into $\alpha'_1,\cdots , \alpha'_\ell$, we have a tiling of $\alpha''_i$ by $\mathcal S$.
%

The tiling rules of $\alpha'_i$s by $\mathcal S$ are called {\it substitution rules}.  A tiling of ${\bf E}^2$ obtained by substitution rules is called {\it a hierarchical tiling}, and if it doesn't have no periodicity of parallel transformation, it is called {\it non-periodic}.
%

An L-tiling is an example of a non-periodic hierarchical tiling. Let $\alpha=$\ \raisebox{-2pt}{\includegraphics[height=12pt]{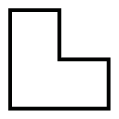}} and $\mathcal{S}=\{ \alpha \}$.  Let $\lambda=2$ and $\alpha'=$\ \raisebox{-4pt}{\includegraphics[height=24pt]{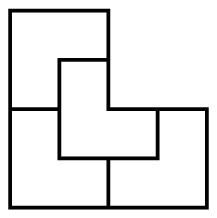}}, then it gives a substitution rule.  
We call this tiling {\it an L-tiling} and this rule {\it L-substitution}.
%

\begin{definition}[$s$-spread]
We call $\alpha'$, a tiling of once L-substitution from $\alpha$, {\it $1$-spread}.  We call $\alpha^{(s)}$, a tiling of $s$ times L-substitution from $\alpha$, {\it $s$-spread}. 
\end{definition}
\par\bigskip
Next, we define a edge and a perturbed edge.  
\begin{definition}[edge]
Let {\it an edge} be a pair of a segment and a one-side neighborhood.  See Figure 3.  
\end{definition}
\begin{figure}
\includegraphics[]{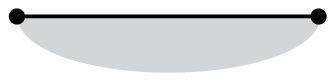}
\caption{an edge}
\end{figure}
For an edge, we regard a segment as a edge of a prototile, and one-side neighborhood as inside of a prototile. 
In the sequel, we do not distinguish an edge of a prototile from an edge in this definition.
Next, we define a perturbed edge. 
\begin{definition}[perturbed edge]
For an edge, fixing the both ends of the edge and perturbing it a little, we get {\it a perturbed edge}.  See Figure 4.   For two perturbed edges $a,b$, $a=b$ if they are congruent.
\end{definition}

\begin{figure}
\includegraphics[]{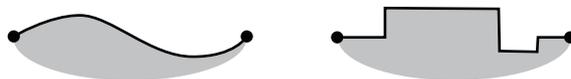}
\caption{examples of perturbed edges}
\end{figure}

In order to perturb edges of prototiles, there exists restriction on a way of perturbing,  because each prototile must be connected.  In our context, we only concern degree of freedom of perturbed prototiles, so we consider only a perturbed edge which is a little perturbed to avoid the restriction.  
%

Next, we define a product of perturbed edges.
\begin{definition}[product of edges]
Let $a_1,a_2,\cdots,a_k$ be perturbed edges.  If they are placed on a straight line from right to left and form a row, then we call it a product of $a_1,a_2,\cdots,a_k$ and we denote this product by $a_1a_2\cdots a_k$.  See Figure 5. 
\end{definition}
\begin{figure}
\includegraphics[]{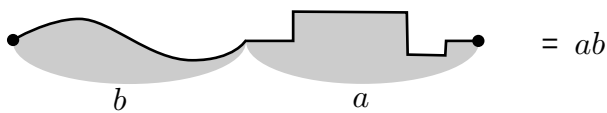}
\caption{an example of a product of perturbed edges}
\end{figure}
For a perturbed edge $a$, we define two operations $\overline{a}$, and $a^{-1}$.  For a perturbed edge $a$, $\overline{a}$ is a symmetry (right-side-left) image of $a$. In the same way, $a^{-1}$ is an upside-down image of $a$.  See Figure 6.
\begin{figure}
\includegraphics[]{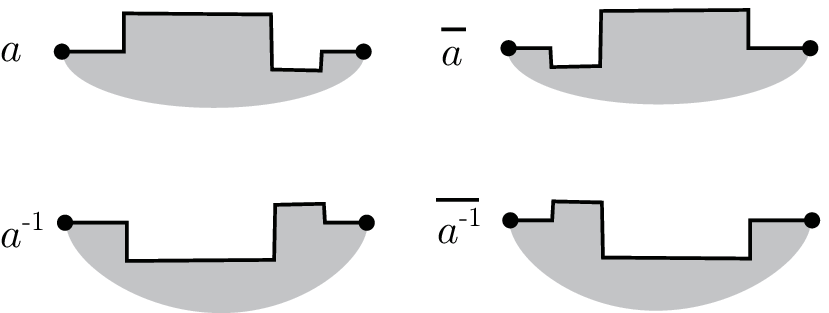}
\caption{definition of $\overline{a}$, and $a^{-1}$ }
\end{figure}
It is easy to show the following lemma.
\begin{lem}
(1) $\overline{\left(\overline{a}\right)}=a, (a^{-1})^{-1}=a$\par
(2) $\overline{ab}=\overline{a} \overline{b}, (ab)^{-1}=b^{-1}a^{-1}$\par
(3) $\overline{(a^{-1})}=\left(\overline{a}\right)^{-1}$
\end{lem}
%

In a tiling, if a tile with a perturbed edge $a$ and another tile with a perturbed edge $b$ are neighbors at $a$ and $b$,  we have $a=\overline{b^{-1}}$.  We denote this relation by $\dfrac{a}{b}$.  We often say that $a$ matches $b$.
%

The following lemma is trivial.
\begin{lem}
(1) $\dfrac{a}{b}$ if and only if $\dfrac{b}{a}$\par
(2) If $\dfrac{a}{b}$ and $\dfrac{a}{c}$ then $b=c$ \par
(3) $\dfrac{ab}{cd}$ if and only if $\dfrac{a}{d}$ and $\dfrac{b}{c}$
\end{lem}
%

Let $\mathcal T$ be a tiling with respect to a tile set $\mathcal S$.  
Suppose that all prototiles are polygons.  Here we assume that there is no vertex of a tile lying on an edge of another tile. 
\begin{definition}[escherization, escher degree]
(1) Let $\mathcal T$ and $\mathcal S$ be as above.  If we perturb edges of prototiles such that the perturbed prototiles give another tiling, we call this process {\it escherization}. \par
(2) If the set of escherization of $\mathcal T$ is parametrized by some perturbed edges, {\it the escher degree} is the number of the parameters.
\end{definition}
%

\begin{example}[escher degree of (P1)]
Let $\alpha$ be a parallelogram and (P1) a tiling of ${\bf E}^2$ as in Figure 2.  Let $a,b,c,d$ be edges of $\alpha$ as in Figure 7.
\begin{figure}
\includegraphics[]{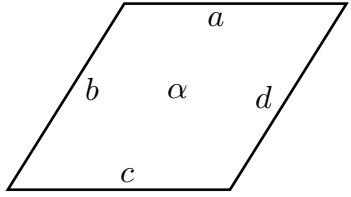}
\caption{Edges of the parallelogram $\alpha$}
\end{figure}

From the matching ofr the tiling, we have $\dfrac{a}{c}$ and $\dfrac{b}{d}$.  That is, if we perturb $a$, then the edge $c$ changes such that $c=\overline{a^{-1}}$, and we can perturb $b$ independently of $a$.  Then the edge $d$ changes such that $d=\overline{b^{-1}}$.  See Figure 8.  
We call relations obtained from the tiling property {\it edge-matchings}.
\begin{figure}
\includegraphics[]{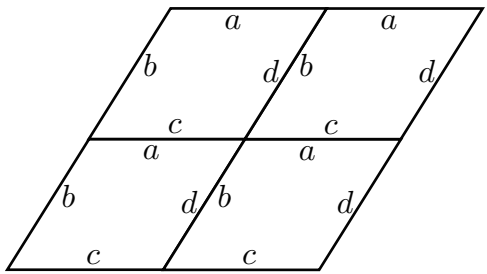}
\caption{matching of the tiling (P1)}
\end{figure}
Hence all escherization of the tiling (P1) is parametrized by edges $a$ and $b$.  So the escher degree is $2$.
\end{example}
%


\section{Escher degree of L-tiling by one prototile}

In this section we show that the escher degree of L-tiling by one prototile is one.  We assume that an L-figure prototile $\alpha$ has $8$ perturbed edges $a,b,c,\cdots, h$ as in Figure 9.  
\begin{figure}
\includegraphics[]{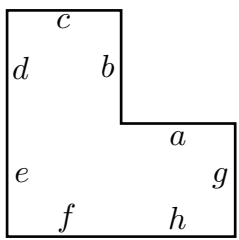}
\caption{prototile $\alpha$}
\end{figure}
%

\begin{thm}[Theorem A]
(1) In a non-periodic L-tiling by one prototile $\alpha$, we have
$$ \frac{a=c=e=g}{b=d=f=h}.$$
(2) The escher degree of this tiling is one. See Figure 18.
\end{thm}
\begin{remark}
$ \frac{a=c=e=g}{b=d=f=h} $ means $a=c=e=g$, $b=d=f=h$, and $\frac{a}{b}$
\end{remark}
%

{\bf Proof:}\quad (1) Considering $1$-spread, we directly have 
$$ \dfrac{b}{c}, \dfrac{a}{d}, \dfrac{h}{c}, \dfrac{b}{e}, \dfrac{f}{a}, \dfrac{h}{c}, \dfrac{g}{b}, \dfrac{h}{a}.$$
 (See Figure 10.) This follows that $ \frac{a=c=e=g}{b=d=f=h}$.
\begin{figure}
\includegraphics[]{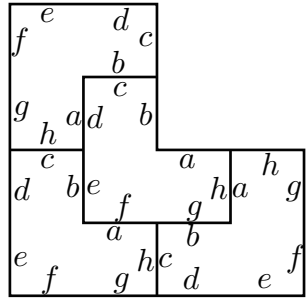}
\caption{matching in $1$-spread of $\alpha'$}
\end{figure}
%

(2) The proof of (1) implies the following lemma.
\begin{lem}\label{lem:3-2}
There exists a $1$-spread of $\alpha$ if and only if $\alpha$ satisfies $\frac{a=c=e=g}{b=d=f=h}$.  
\end{lem}
Let $\alpha'$ be a $1$-spread and $a',b',\cdots, h'$ its edges. (See Figure 11.) 
\begin{figure}
\includegraphics[]{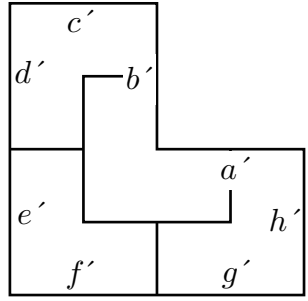}
\caption{perturbed edges of $\alpha'$}
\end{figure}

Then we have 
$$ a'=ha, b'=bc,c'=e'=g'=de, d'=f'=h'=fg.$$
(See Figure 10.)  If $\frac{a=c=e=g}{b=d=f=h}$ then we easily show that 
$$ \frac{a'=c'=e'=g'}{b'=d'=f'=h'} .$$
(For example, $h=d$ and $a=e$ implies $a'=ha=de=c'$, $\frac{h}{c}$ and $\frac{a}{b}$ implies $\frac{ha}{bc}$ and $\frac{a'}{b'}$.)  This follows that a $1$-spread of $\alpha'$ exists, that is, a $2$-spread of $\alpha$ exists.

In the same way, if we assume that we have an $s$-spread $\alpha^{(s)}$ of $\alpha$ and $\alpha^{(s-1)}$ has edges $a^{(s-1)}$, $b^{(s-1)}$, $\cdots$, $h^{(s-1)}$, then it satisfies that $$\dfrac{a^{(s-1)}=c^{(s-1)}=e^{(s-1)}=g^{(s-1)}}{b^{(s-1)}=d^{(s-1)}=f^{(s-1)}=h^{(s-1)}}.$$
If we set 
\begin{gather*} a^{(s)}=h^{(s-1)}a^{(s-1)}, b^{(s)}=b^{(s-1)}c^{(s-1)},\\
c^{(s)}=e^{(s)}=g^{(s)}=d^{(s-1)}e^{(s-1)}, \\
d^{(s)}=f^{(s)}=h^{(s)}=f^{(s-1)}g^{(s-1)}\end{gather*}
inductively, then they satisfy 
$$ \frac{a^{(s)}=c^{(s)}=e^{(s)}=g^{(s)}}{b^{(s)}=d^{(s)}=f^{(s)}=h^{(s)}}.$$
and it follows that an $1$-spread of $\alpha^{(s)}$ exists, that is, we have an $(s+1)$-spread of $\alpha$.

If $\alpha$ satisfies $\dfrac{a=c=e=g}{b=d=f=h}$ then an $s$-spread $\alpha^{(s)}$ exists for any $s$.  Hence the escherization is parametrized by $a$ and the escher degree is one.  This completes the proof. 
%


\section{L-tiling by two prototiles (1)}

In this section, we consider the cases where two prototiles $\alpha,\beta$ (Figure 12) make one $1$-spread $\alpha'$. 
\begin{figure}
\includegraphics[]{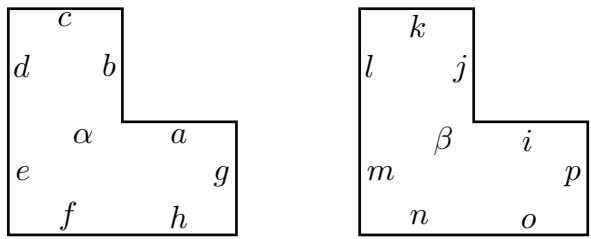}
\caption{two prototiles $\alpha$ and $\beta$}
\end{figure}

From Theorem A, if $\alpha'$ satisfies $\frac{a'=c'=e'=g'}{b'=d'=f'=h'}$ then 
$\alpha'$ has $s$-spread for $s=1,2,3,\cdots$.  We observe 5 patterns of $\alpha'$ in Figure 13.  We call them no.8, no.4, no.2, no.3, and no.5 respectively.  (The numbering order is not ascending nor descending.  These numberings are determined by the order of $\alpha$ and $\beta$.)  
\begin{figure}
\includegraphics[]{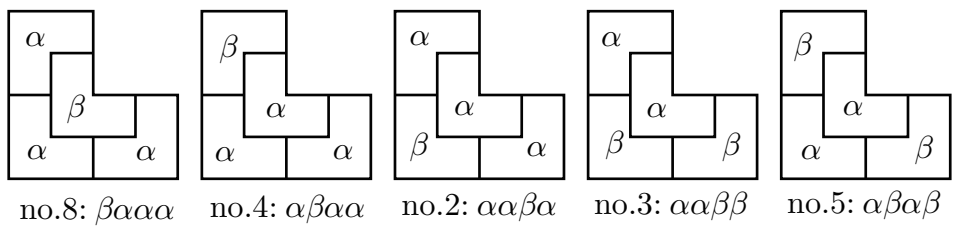}
\caption{five patterns of $1$-spread}
\end{figure}

And we have the following theorem.
\begin{thm}[Theorem B]
(1) For no.2, no.3, no.4, the escher degree is $1$.  \par
(2) For no.5, we have $$ \frac{a=e=m}{b=f=n}, \frac{c=g=i=k=o}{d=h=j=l=p} $$ and the escher degree is $2$. (See Figure 19.)\par
(3) For no.8, we have $$ \dfrac{a}{l=n=p},\dfrac{k=m=o}{b},\dfrac{c=g}{d=h},\dfrac{e=i}{f=j}$$ and the escher degree is $4$.  (See Figure 20.)
\end{thm}

{\bf Proof:} \quad If $\alpha'$ satisfies $\frac{a'=c'=e'=g'}{b'=d'=f'=h'}$, then $\alpha'$ has $s$-spread for any $s$. So, it is sufficient to solve  the edge-matching in $\alpha'$ and $\frac{a'=c'=e'=g'}{b'=d'=f'=h'}$.
\begin{figure}
\includegraphics[]{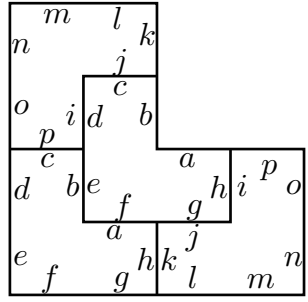}
\caption{$\alpha'$ for no.5 tiling}
\end{figure}
For example, for no.5, edge-matching is given as in Figure 14 and we have
$$  \dfrac{j}{c},\frac{i}{d},\frac{p}{c},\frac{b}{e},\frac{f}{a},\frac{h}{k},\frac{g}{j},\frac{h}{i}.$$
$\frac{a'=c'=e'=g'}{b'=d'=f'=h'}$ implies 
$$ \frac{pa=lm=de=lm}{bk=no=fg=no} .$$
We solve the system of equation and we have $ \frac{a=e=m}{b=f=n}, \frac{c=g=i=k=o}{d=h=j=l=p} $. Inversely if $ \frac{a=e=m}{b=f=n}, \frac{c=g=i=k=o}{d=h=j=l=p} $ then $\frac{a'=c'=e'=g'}{b'=d'=f'=h'}$ is satisfied.  For other tilings, we can solve the system of relations in a similar way.

\par\bigskip
\begin{remark}
If we  have 
$$ \frac{a=b=e=g=i=k=m=o}{b=d=f=h=j=l=n=p}, $$
it follows that $\alpha=\beta$.  So we can return the case of one prototile $\alpha$ and the escher degree is 1.  This means that there are no solution of two distinct prototiles.  In the two prototiles case, if the escher degree is more than $1$, we call the tiling non-trivial. There are two types (no.5 and no.8) of non-trivial tilings by two prototiles and one $1$-spread.
\end{remark}

\section{L-tiling by two prototiles (2)}

In this section, we consider two prototiles $\alpha,\beta$ (Figure 12) and make two $1$-spreads $\alpha', \beta'$.  There are 16 possibilities for $\alpha'$ and so as for $\beta'$. 

We determine numbering for $\alpha'$ and $\beta'$ as in Figure 15, and we represent a tiling by a pair of two numbers for $\alpha'$ and $\beta'$.  For example, $(5,10)$ means that $\alpha'$ and $\beta'$ given in Figure 16.
\begin{figure}
\includegraphics[]{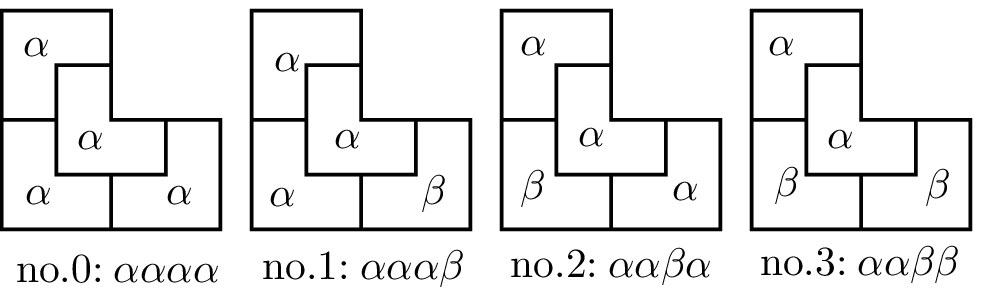}
\includegraphics[]{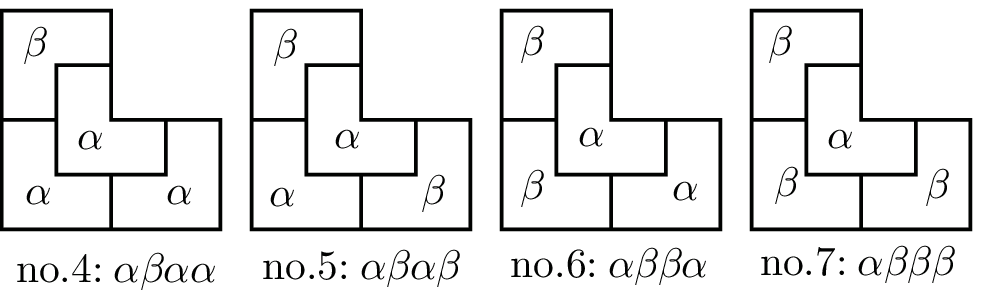}
\includegraphics[]{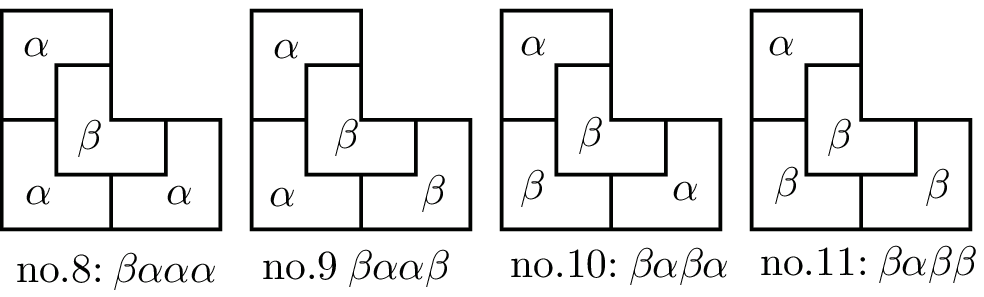}
\includegraphics[]{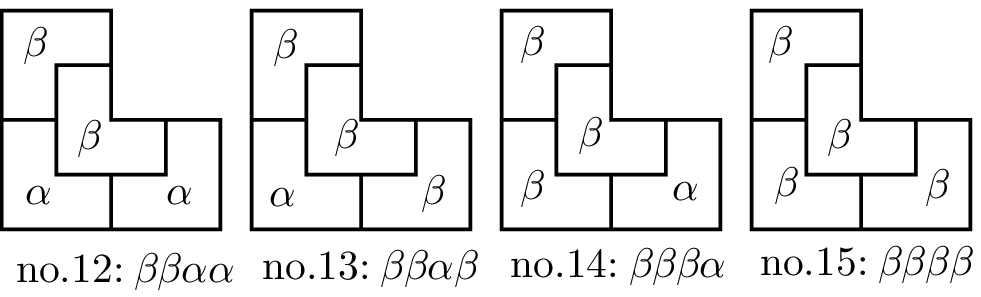}
\caption{numbering for $\alpha'$ and $\beta'$}
\end{figure}
\begin{figure}
\includegraphics[]{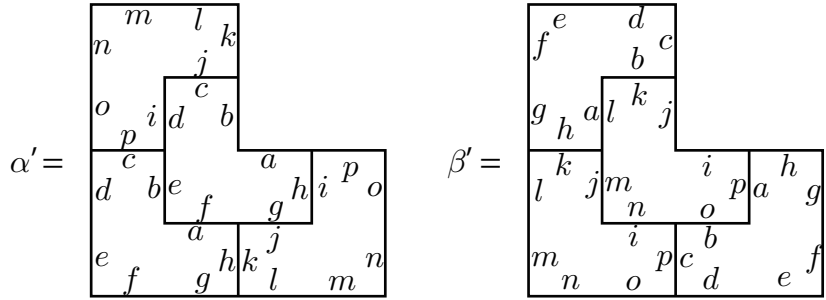}
\caption{ Tiling $(5,10)$}
\end{figure}
We remove the case $\alpha'=\beta'$ and two trivial cases $(0,15)$ and $(15,0)$, there remains $238$ combinations. If we exchange the role of $\alpha$ and $\beta$, we know that $(i,j)$ and $(15-j,15-i)$ are equivalent.  From the following lemma, we conclude that the number of remaining combinations is $119$. 

\begin{lem}
$(i,15-i)$ and $(15-i,i)$ are equivalent.
\end{lem}
{\bf Proof:}\quad If we denote $\alpha_{(i,j)}^{(s)}$ (resp. $\beta_{(i,j)}^{(s)}$) 
by the $s$-spread of $\alpha$ (resp. $\beta$) of tiling $(i,j)$, it is easily show that 
\begin{align*}
&\alpha_{(i,15-i)}^{(2k-1)}=\beta_{(15-i,i)}^{(2k-1)},
&\beta_{(i,15-i)}^{(2k-1)}=\alpha_{(15-i,i)}^{(2k-1)}\\
&\alpha_{(i,15-i)}^{(2k)}=\alpha_{(15-i,i)}^{(2k)},
&\beta_{(i,15-i)}^{(2k)}=\beta_{(15-i,i)}^{(2k)},
\end{align*}
for any $k=1,2,\cdots$. This completes the proof. \par\bigskip
Here we have the third theorem.
\begin{thm}[Theorem C]
(1) If the escher degree of $(i,j)$ tiling is more than $1$, then $(i,j)=(5,10),(10,5),(0,2),(13,15),(0,8),(7,15)$, $(0,10)$, $(5,15)$.\par
(2) For the tiling $(5,10)$ (equivalently $(10,5)$), $\frac{a=c=g=m}{f=j=l=p},\frac{e=i=k=o}{b=d=h=n}$ and the escher degree is $2$.  (See Figure 21.)\par
(3) For the tiling $(0,2)$ (equivalently $(13,15)$), $\frac{a=c=e=g=i=k}{b=d=f=h=j=p}$ and the escher degree is $5$.  (See Figure 22.)\par
(4) For the tiling $(0,8)$ (equivalently $(7,15)$), $\frac{a=c=e=g=k=m=o}{b=d=f=h=l=n=p}$ and the escher degree is $3$.  (See Figure 23.)\par
(5) For the tiling $(0,10)$ (equivalently $(5,15)$), $\frac{a=c=e=g=k=o}{b=d=f=h=l=p},\frac{m}{j},\frac{i}{n}$ and the escher degree is $3$.  (See Figure 24.)
\end{thm}

For each $(i,j)$, we solve a system of equations of edge-matchings of $s$-spread ($s=1,2,\cdots$).   

In some cases, only $i$ (resp. only $j$) determines the result.  For example, the following lemma holds.
\begin{lem}
The escher degree of $(1,j)$ is $1$ for any $j$.
\end{lem}
\begin{figure}
\includegraphics[]{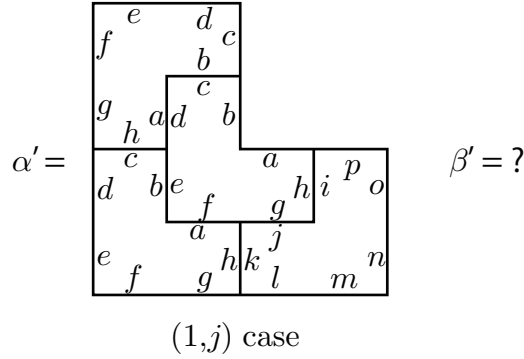}
\caption{$(1,j)$ tiling}
\end{figure}
{\bf Proof:}  Assume that $\alpha'$ is no.1.  (See Figure 17.)  From the edge-matching of $\alpha'$, we get 
$$ \dfrac{c=e=k=i}{b=h},\dfrac{a}{d=f},\dfrac{g}{j}. $$
Edges of $\alpha'$ is given by $a'=pa,b'=bc,c'=e'=de,d'=f'=fg,g'=lm,h'=no$ and we have 
$$ \dfrac{c'=e'}{b'=h'},\dfrac{a'}{d'=f'}, $$ 
in the edge-matching of $\alpha''$, and we get additional conditions $\dfrac{d}{c}, b=n, c=o, \dfrac{p}{g}$ and hence $\dfrac{a=c=e=k=i=o}{b=d=f=h=n},\dfrac{g}{j=p}$.  In the edge-matching in $\alpha'''$, we obtain another condition $\dfrac{d'}{c'}$ , hence $\dfrac{d}{g},\dfrac{e}{f}$, and we have \\
$\dfrac{a=c=e=g=k=i=o}{b=d=f=h=j=n=p}$.  In the edge-matching in $\alpha''''$, we have $\dfrac{d'}{g'}$, hence $\dfrac{f}{m}$ and $\dfrac{g}{l}$, and we have $\dfrac{a=c=e=g=i=k=m=o}{b=d=f=h=j=l=n=p}.$ This implies that within at most $4$-spread all edges $a,b,\cdots, p$ are related and $\alpha=\beta$.  This completes the proof. 
\par\bigskip
And in a similar way, we can show that the escher degree of any tiling other than $(i,j)=(5,10),(10,5),(0,2),(13,15),(0,8),(7,15)$, $(0,10)$, $(5,15)$ is $1$.

(2) In $(i,j)=(5,10)$ case, the edge-matching of $\alpha'$ is
\begin{equation} \frac{j}{c}, \frac{i}{d},\frac{p}{c},\frac{b}{e},\frac{f}{a},\frac{h}{k},\frac{g}{j},\frac{h}{i}. \tag{a}
\end{equation}
(See the left figure of Figure 16.)  The  edge-matching of $\beta'$ is
\begin{equation} \frac{b}{k},\frac{a}{l},\frac{h}{k},\frac{j}{m},\frac{n}{i},\frac{p}{c},\frac{o}{b},\frac{p}{a}. 
\tag{b}
\end{equation}
(See the right figure of Figure 16.) (a),(b) are equivalent to
\begin{equation}
\frac{a=c=g=m}{f=j=l=p},\frac{e=i=k=o}{b=d=h=n}.\tag{c}
\end{equation}
Let $a',b',\cdots, o',p'$ be edges of $\alpha',\beta'$ as in Figure 11, we have
\begin{align*}
&a'=pa, b'=bk,c'=lm,d'=no,e'=de,f'=fg,g'=lm,h'=no\\
&i'=hi,j'=jc,k'=de,l'=fg,m'=lm,n'=no,o'=de,p'=fg
\end{align*}
Since the formula (c) holds for $\alpha', \beta'$, we have
$$ \frac{a'=c'=g'=m'}{f'=j'=l'=p'},\frac{e'=i'=k'=o'}{b'=d'=h'=n'}.$$ 
For example, $p=l$ and $a=m$ implies $a'=pa=lm=c'$, and so on.  These are edge-matching of $2$-spread $\alpha'',\beta''$. 

Inductively, we observe as follows.  
Suppose that we have $\alpha^{(s)}, \beta^{(s)}$.  Let 
$a^{(s)},b^{(s)},\cdots,o^{(s)}, p^{(s)}$ be edges of $\alpha^{(s)}, \beta^{(s)}$.  The relation between $a^{(s-1)}$s and $a^{(s)}$s are given by 
\begin{align*}
&a^{(s)}=p^{(s-1)}a^{(s-1)}, b^{(s)}=b^{(s-1)}k^{(s-1)},
c^{(s)}=l^{(s-1)}m^{(s-1)},d^{(s)}=n^{(s-1)}o^{(s-1)},\\
&e^{(s)}=d^{(s-1)}e^{(s-1)},f^{(s)}=f^{(s-1)}g^{(s-1)},
g^{(s)}=l^{(s-1)}m^{(s-1)},h^{(s)}=n^{(s-1)}o^{(s-1)}\\
&i^{(s)}=h^{(s-1)}i^{(s-1)},j^{(s)}=j^{(s-1)}c^{(s-1)},
k^{(s)}=d^{(s-1)}e^{(s-1)},l^{(s)}=f^{(s-1)}g^{(s-1)},\\
&m^{(s)}=l^{(s-1)}m^{(s-1)},n^{(s)}=n^{(s-1)}o^{(s-1)},
o^{(s)}=d^{(s-1)}e^{(s-1)},p^{(s)}=f^{(s-1)}g^{(s-1)},
\end{align*}
for any $s=0,1,2,\cdots$.  Using simple calculations, we have the following lemma.

\begin{lem} \mbox{}\vspace{3pt}\par
If 
$\dfrac{a^{(s-1)}=c^{(s-1)}=g^{(s-1)}=m^{(s-1)}}{f^{(s-1)}=j^{(s-1)}=l^{(s-1)}=p^{(s-1)}},\dfrac{e^{(s-1)}=i^{(s-1)}=k^{(s-1)}=o^{(s-1)}}{b^{(s-1)}=d^{(s-1)}=h^{(s-1)}=n^{(s-1)}}$, then\\ $\dfrac{a^{(s)}=c^{(s)}=g^{(s)}=m^{(s)}}{f^{(s)}=j^{(s)}=l^{(s)}=p^{(s)}},\dfrac{e^{(s)}=i^{(s)}=k^{(s)}=o^{(s)}}{b^{(s)}=d^{(s)}=h^{(s)}=n^{(s)}}$.
\end{lem}

{\bf Proof:}\quad For example, $p^{(s-1)}=l^{(s-1)}$ and $a^{(s-1)}=m^{(s-1)}$ implies $a^{(s)}=p^{(s-1)}a^{(s-1)}=l^{(s-1)}m^{(s-1)}=c^{(s)}$.  Other relations are shown in a similar way.
From this lemma, $(s+1)$-spreads $\alpha^{(s+1)},\beta^{(s+1)}$ exist for any $s$. 

(3), (4), and (5) are shown in a similar way as (2). 

\begin{remark}
In Appendix, we show figures of these tilings.  In $(0,2)$, $(0,8)$ tilings, $\beta^{(s)}$ contains only one tile $\beta$.   This means that these tilings are equivalent to a tiling of $\alpha$ as infinite tilings.  
\end{remark}


\section{Appendix (1)}

From Figure 18 to Figure 24 are pictures of tilings appearing in Theorems A,B, and C.

\section{Appendix (2)}

For a tile set $\mathcal T$, if {\bf any} tiling by $\mathcal T$ has no periodicity of parallel transformation then we call $\mathcal T$ {\it an aperiodic tile set.}  Any tilings we obtain in this paper are not aperiodic tile sets.  From Figure 25 to Figure 28 are figures of periodic tilings.

\begin{figure}[]
\includegraphics[]{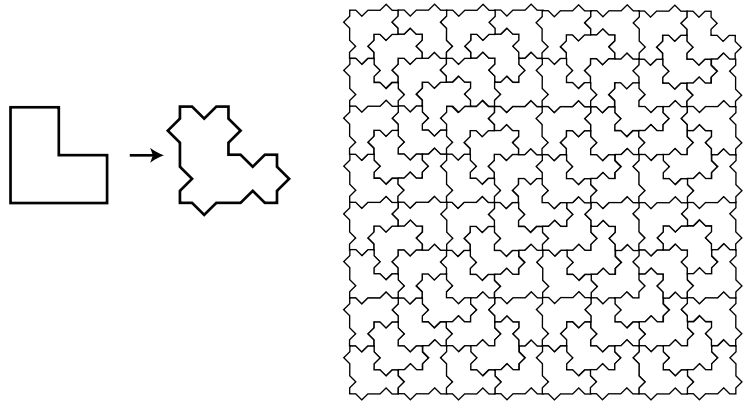}
\caption{Tiling of one prototile.}
\end{figure}

\begin{figure}[]
\includegraphics[]{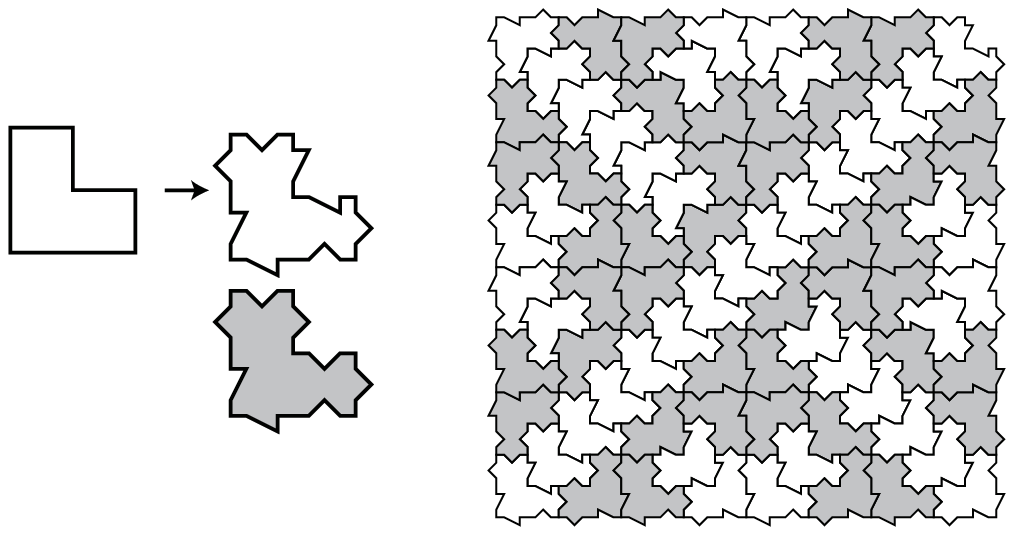}
\caption{Tiling of two prototiles, no.5.}
\end{figure}

\begin{figure}[]
\includegraphics[]{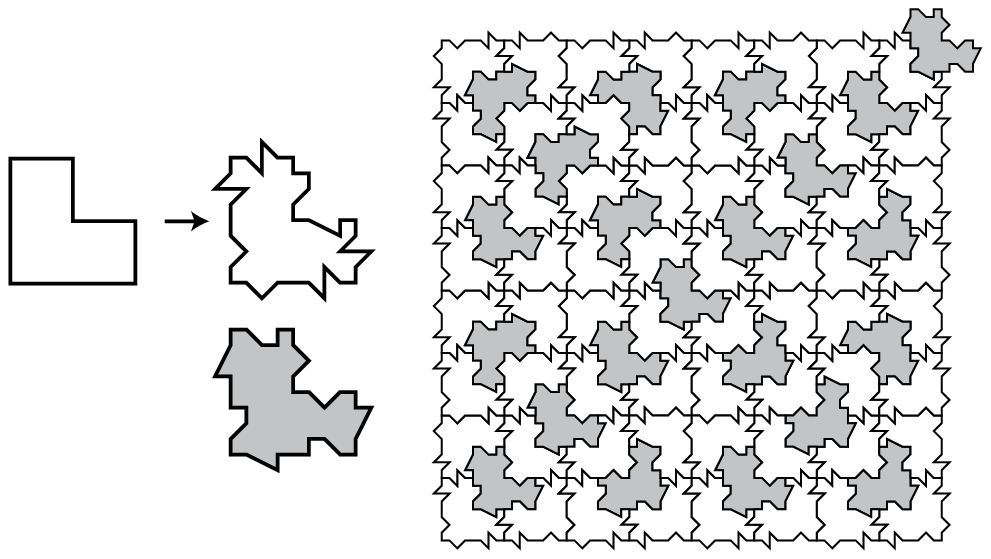}
\caption{Tiling of two prototiles, no.8.}
\end{figure}

\begin{figure}[]
\includegraphics[]{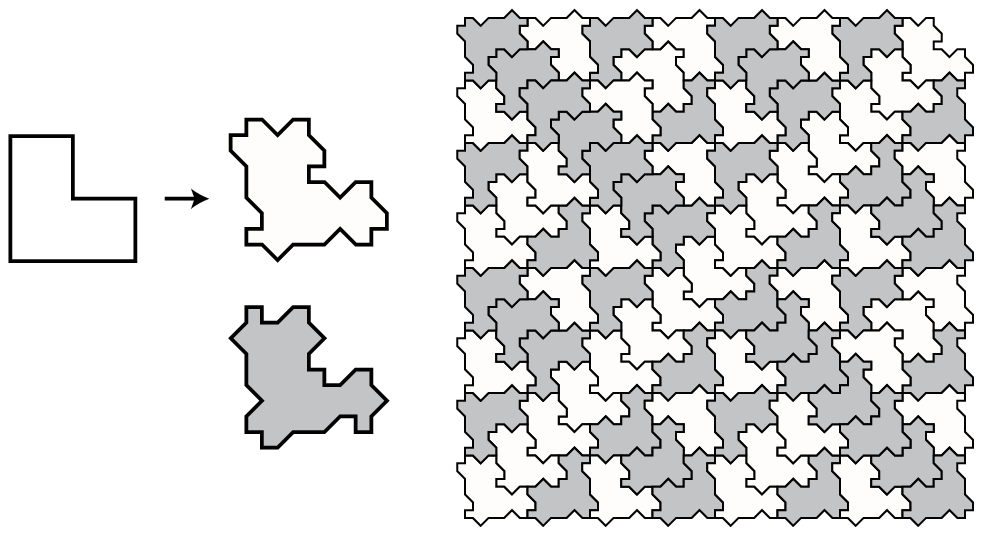}
\caption{Tiling of two prototiles (5,10).}
\end{figure}

\begin{figure}[]
\includegraphics[]{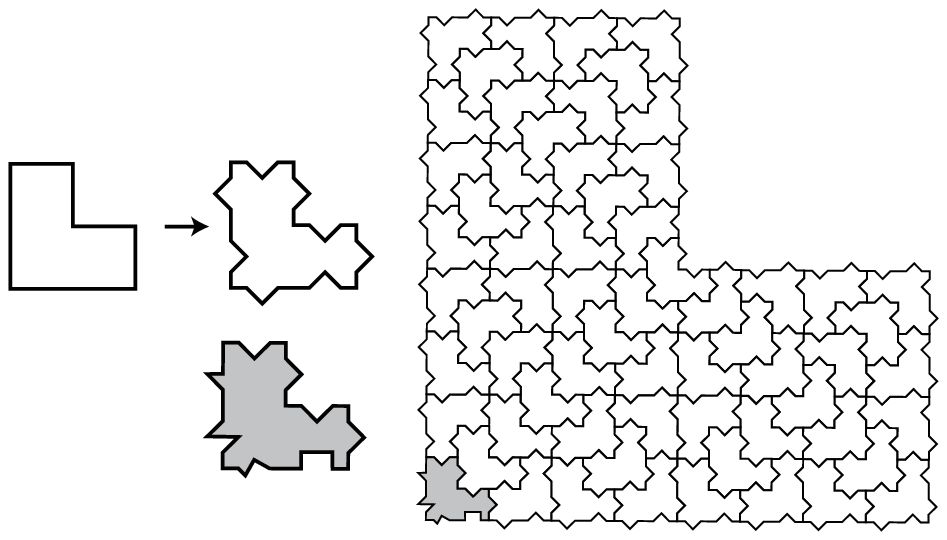}
\caption{Tiling of two prototiles (0,2).}
\end{figure}

\begin{figure}[]
\includegraphics[]{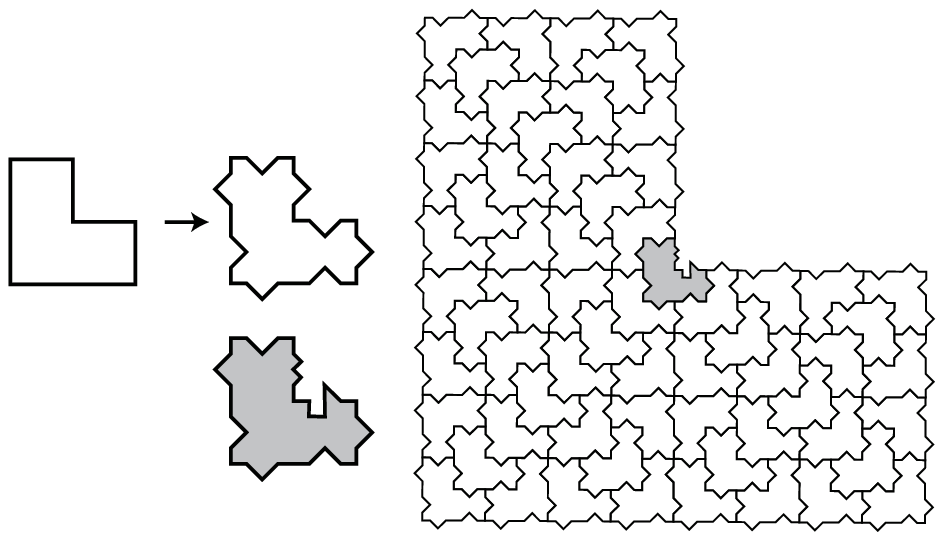}
\caption{Tiling of two prototiles (0,8).}
\end{figure}

\begin{figure}[]
\includegraphics[]{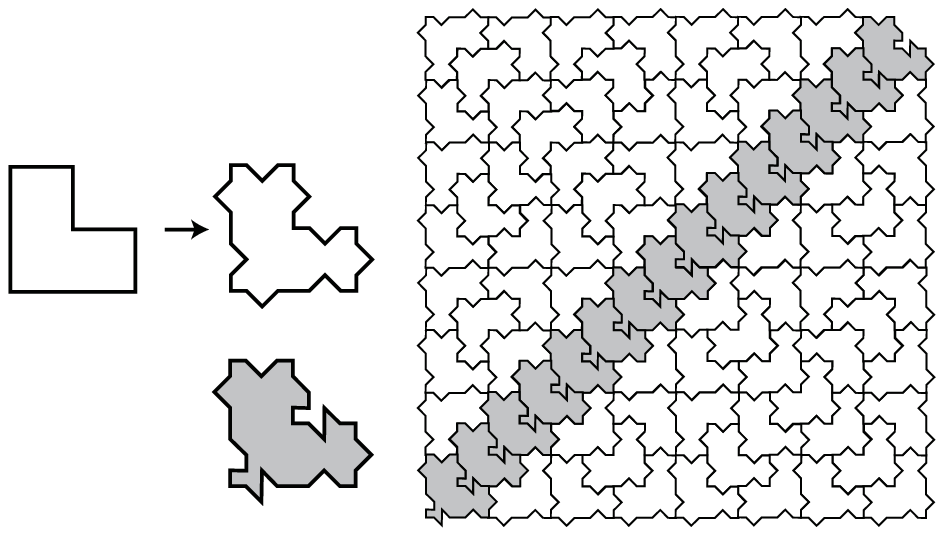}
\caption{Tiling of two prototiles (0,10).}
\end{figure}

\begin{figure}[]
\includegraphics[]{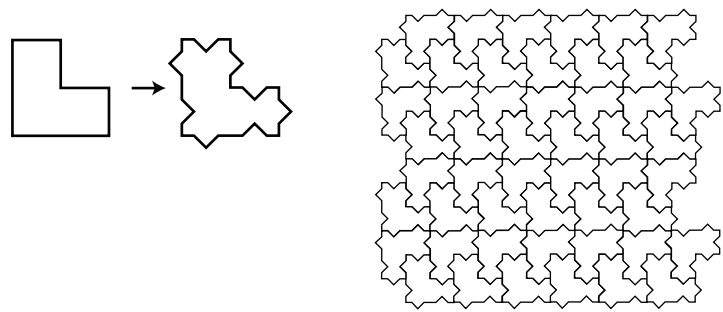}
\caption{Periodic tiling of one prototile.}
\end{figure}

\begin{figure}[]
\includegraphics[]{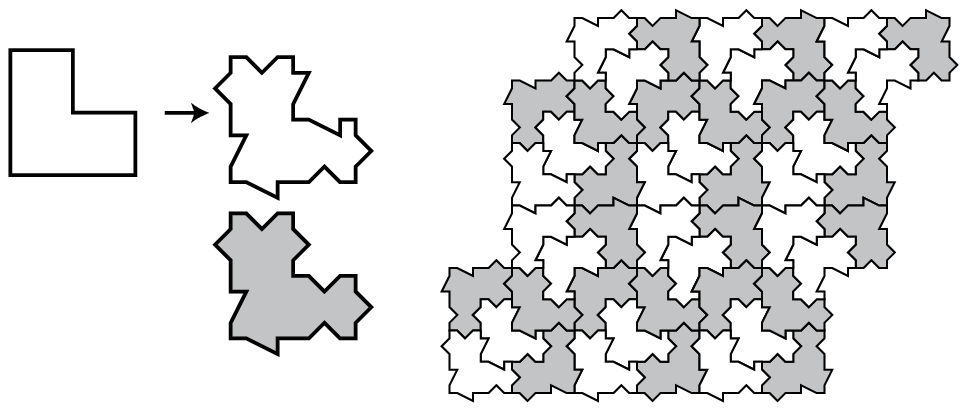}
\caption{Periodic tiling of two prototiles no.5.}
\end{figure}

\begin{figure}[]
\includegraphics[]{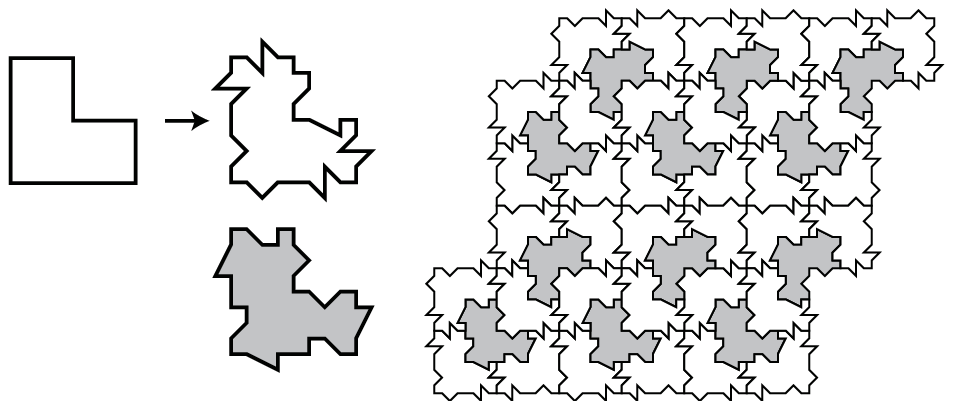}
\caption{Periodic tiling of two prototiles no.8.}
\end{figure}

\begin{figure}[]
\includegraphics[]{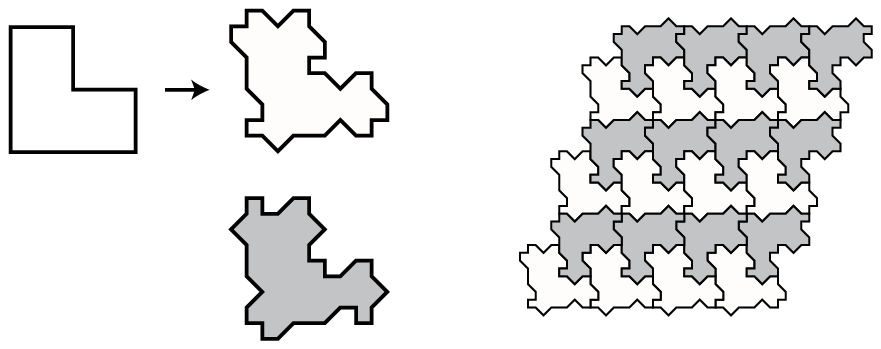}
\caption{Periodic tiling of two prototiles (5,10).}
\end{figure}

\bibliographystyle{amsplain}

\end{document}